\documentclass[11pt]{article}
\usepackage{amsmath}
\usepackage{amsfonts}
\usepackage{amssymb}

\newtheorem{theorem}{Theorem}

\newcommand{\R}{\mathbb{R}}
\newcommand{\C}{\mathbb{C}}
\newcommand{\Z}{\mathbb{Z}}
\newcommand{\CP}{\mathbb{CP}}
\newcommand{\M}{M_{KT}}
\newcommand{\KT}{Kodaira--Thurston }

\begin{document}
\title {Theta functions on the Kodaira--Thurston manifold
\thanks{This work is supported by RFFI, grant 06-01-00094-a.}}
\author{D.V. Egorov}
\date{}
\maketitle

\begin{abstract}
We define  analogue of theta functions on the Kodaira--Thurston manifold which is a compact 4-dimensional symplectic manifold and use them to construct canonical symplectic embedding of the Kodaira--Thurston manifold into the complex projective space (analogue of the Lefshetz theorem).
\end{abstract}


\section{Introduction}
In this work we construct analogue of classical theta functions on the
Abelian variety for the Kodaira--Thurston nilmanifold $\M$. The
classical theta function from the geometrical point of view is a
section of  a holomorphic line bundle over complex torus. The theta
function on the Kodaira--Thurston manifold is defined as section (not
holomorphic) of a special line bundle $L$ over $\M$.

Analogues of theta functions for nilmanifolds were defined earlier
\cite{A,KU}, but all these generalizations are based on
representation theory. We construct such {\sl analogue of  theta
functions with characteristics that they set canonical symplectic
embedding of $\M$ into a complex projective space (analogue of
Lefshetz theorem)}.

The Kodaira--Thurston manifold $\M$ is the quotient of $\R^4$ by the
free action of discrete group $\Gamma$, which generators are
\begin{equation}
\label{period}
\begin{array}{clll}
a: (x,y,z,t) & \to & (x+1, y, z+y, t)\\
b: (x,y,z,t) & \to & (x, y+1, z, t) \\
c: (x,y,z,t) & \to & (x, y, z+1, t) \\
d: (x,y,z,t) & \to & (x, y, z, t+1)
\end{array}
\end{equation}
The \KT manifold is notable as the first known example of a
symplectic but not K\"{a}hler manifold \cite{Thurston}.

Note that the embedding of $\M$ into $\CP^n$ can't be holomorphic,
since $\M$ is not K\"{a}hler. However we prove that this map is
symplectic. In other words the Fubini--Study form on $\CP^n$ induces
symplectic structure on $\M$.

In \S 2 we recall necessary facts from classical theta function
theory, in \S 3 we give the definition of theta function on $\M$ and
study some of its properties, in \S 4 we construct embedding into
complex projective space (Theorem 1) and in \S 5 prove that this
embedding is symplectic (Theorem 2).

Author is grateful to I.A. Taimanov for posing the problem and to
A.E. Mironov for helpful discussions.

\section{Classical theta function}
Let's recall some known facts about the classical theta function of Jacobi on one-dimensional complex torus.
We will need them in the sequel.

Consider formal series
$$
\theta(z,\tau) =
\sum_{k\in\mathbb{Z}}{e^{2\pi i kz + \pi ik(k-1)\tau }} .
$$
If $\mathrm{Im}\,\tau >0$ then this series converges uniformly on every compact subset of $\C$ and thus define entire function. In notations of \cite{Mumford} this theta function is written as
$$
\exp(-\pi i\tau/4 + \pi i z)\cdot\theta_{-1/2,0}(z,\tau).
$$

The invariance under map $\tau\to \tau + 1$ and certain periodicity conditions explain our particular choice
\begin{eqnarray}
\label{classperiod}
\theta(z+1, \tau)& = &\theta(z, \tau), \\
\label{classperiod2}
\theta(z+\tau, \tau)& = & \exp({-2\pi iz })\cdot \theta(z,\tau).
\end{eqnarray}

The generalization of theta function is the theta function of degree $k\in\mathbb{N}$. It is an entire function $\theta_k(z,\tau)$ which periodicity conditions are as follows
$$
\theta_k(z+1, \tau) =  \theta_k(z, \tau),
$$
$$
\theta_k(z+\tau, \tau) =   \exp({-2\pi i kz})\theta_k(z,\tau).
$$
It is not hard to prove that theta functions of degree $k$ form linear space of dimension $k$. Let's denote this space by $\mathcal{L}_k$.

The product of theta functions can be the theta function of higher degree. Let $\{\alpha_i\}^k_{i=1}$ be the set of constants such that the sum of them is equal to zero. Then the following product is theta function of degree $k$:
$$
\prod_{i=1}^k{\theta(z+\alpha_i,\tau)} \in \mathcal{L}_k.
$$

Theta function is equal to zero at the point $z=1/2$. There is a single up to multiplicity zero  in the fundamental region of lattice $\mathbb{Z}+\tau\mathbb{Z}$.

Theta function $\theta(z,\tau)$ satisfies following PDE:
\begin{equation}
\label{classheat}
\frac{\partial\theta}{\partial\tau} = \frac{1}{4\pi i}\cdot\frac{\partial^2\theta}{\partial z^2} - \frac{\theta}{2}
\end{equation}

Let $\{\theta^p_k(z,\tau)\}_{p=1}^k$ be the basis in the space of theta functions of degree $k$. Then map written in homogenous coordinates
$$
\varphi_k(z) = [\theta^1_k(z,\tau):\ \ldots\ :\theta^k_k(z,\tau)]
$$
is well-defined map of complex torus into $\mathbb{CP}^{k-1}$.

The Lefshetz theorem holds: {\sl if $k\geq 3$ then map $\varphi_k$ is embedding}. Note that this theorem is true for Abelian tori of arbitrary dimensions.

\section{Definition of theta function on the \KT manifold}

Projection map $(x,y,z,t) \rightarrow (y,t)$ gives $M_{KT}$ the structure of
$T^2$-bundle over $T^2$. Left-invariant symplectic form $\omega_{{KT}} =
(dz-xdy)\wedge dx + dy\wedge dt$ tames bundle structure. This means that restriction of $\omega_{KT}$ on any fibre and base is not degenerated.

Hence let the space of theta functions of degree $k$ on $\M$ be the span of  products of classical theta functions on  fibre and base:
$$
\theta^p_k(z+ix, y+i)\cdot \theta^{q}_k(y+it, i);\quad p,q = 1,\ldots ,k
$$
Let's denote this space by $\mathcal{L}_k$. Clearly dimension of $\mathcal{L}_k$ is equal to $k^2$.

Theta function of degree one we denote as
$$
\theta_{KT}(x,y,z,t) = \theta(z+ix,y+i)\cdot \theta(y+it,i).
$$

\subsection{Theta function~--- section of complex line bundle}
The periodicity conditions of theta function on $\M$ are as follows:
\begin{equation}
\label{mults}
\begin{array}{lcl}
\theta_{{KT}}(x+1, y, z+y, t) &=& \exp(-2\pi i(z+ix) )\cdot \theta_{{KT}}(x, y, z, t)\\
\theta_{{KT}}(x, y+1, z, t)  &=&   \theta_{{KT}}(x, y, z, t)\\
\theta_{{KT}}(x, y, z+1, t)&=&  \theta_{{KT}}(x, y, z, t)\\
\theta_{{KT}}(x, y, z, t+1) &=& \exp(-2\pi i (y+it) )\cdot \theta_{{KT}}(x, y, z, t)
\end{array}
\end{equation}
These formulae imply that $\theta_{KT}$ is a section of line bundle over $\M$. This bundle is obtained by factorization $\mathbb{R}^4\times \mathbb{C}$ under the action of group $\Gamma$
$$
(u,w) \sim (\lambda\cdot u, e_{\lambda}(u)w),\quad u\in\mathbb{R}^4,w\in\mathbb{C},\lambda \in\Gamma
$$
where $e_\lambda(u)$ are multiplicators i.e. nonzero functions
$$
e_\lambda : \mathbb{R}^4 \rightarrow \mathbb{C}^{*}
$$
which satisfy following identities
$$
e_{\lambda}(\mu\cdot u) e_{\mu}(u) = e_{\lambda\mu}(u),\quad \lambda,\mu\in \Gamma
$$
$$
e_{0}(u) = 1.
$$
Sections of bundle which is set by multiplicators are in one-to-one correspondence with functions $f$ on $\mathbb{R}^4$ such that
$$
f (\lambda \cdot u) = e_\lambda(u) f(u),\quad \lambda\in \Gamma, u\in\mathbb{R}^4.
$$

Let's note that  the relations in group $\Gamma$ hold for multiplicators. It is obvious since multiplicators are set by  behavior of the same function.

\subsection{Multiplicative property of $\theta_{{KT}}$}
Let's define the action of $\zeta=(\zeta^1,\zeta^2)\in\C^2$ on $\theta_{KT}$:
\begin{equation}
\label{preproduct}
(\zeta\cdot \theta_{{KT}})(x,y,z,t) = \theta(z+ix + \zeta^1,y+i)\theta(y+it + \zeta^2,i).
\end{equation}
If
$$
\sum_{i=1}^k{\alpha_i} = 0
$$
then the following product is theta function of degree $k$:
\begin{equation}
\label{product}
\prod_{i=1}^k{(\zeta_i\cdot\theta_{{KT}})(x,y,z,t)} \in \mathcal{L}_k.
\end{equation}
The proof follows from the analogous property of classical theta function.

\section{Embedding of $\M$ into complex projective space}
Let's enumerate the basis of $\mathcal{L}_k$: $\{s_i\}_{i=1}^{k^2}$. Map
$$
\varphi_k = \left(s_1 , s_2 , \ldots , s_{k^2}\right)
$$
is well defined map of $\M$ to $\CP^{k^2-1}$.

\begin{theorem}
If $k\geq 3$ then map $\varphi_k$ is embedding.
\end{theorem}

{\sc Proof.}\quad For brevity we will prove the theorem in case of $k=3$. The proof for $k>3$ is analogous.

Firstly we will prove the injectivity of $\varphi_k$. We will follow the proof of the classical Kodaira embedding theorem stated in \cite[ch. 1, \S
4]{GH}

Note that theta functions of degree $k$ are the global sections of $k$-th tensor power of bundle set by multiplicators \eqref{mults}.

If for any two points $u\neq v\in {M_{KT}}$ there exists section $s\in \mathcal{L}_k$ such that $s(u)=0$ and $s(v)\neq 0$ then map $\varphi_k$ is injective. Indeed, assume that map "glues" points $u$ and $v$. It means that for all sections $s\in\mathcal{L}_k$ it is true that $s(v) = \zeta \cdot s(u)$, where $\zeta$ is nonzero constant. If $s$ satisfies previously mentioned condition  then $\zeta$ must be zero and we have contradiction.

Note also that given condition implies that for any point $u\in{M_{KT}}$ not all sections vanish at point $u$.

We construct needed theta function of degree $k=3$ as a product of two functions $s = fg$:
\begin{multline}
\label{section} f(x,y,z,\alpha,\beta)= \theta(z+ix +\alpha,
y+i)\theta(z+ix +\beta, y+i)\times\\ \times\theta(z+ix -\alpha-\beta, y+i),
\end{multline}
\begin{equation}
\label{sectionend} g(y,t,\gamma,\delta) =  \theta(y+it
+\gamma,i)\theta(y+it+\delta,i)\theta(y+it-\gamma-\delta,i).
\end{equation}
As follows from \eqref{product} product $fg$ is a valid theta function of degree $k=3$ on $\M$.

We denote coordinates of $u,v$ as $(x,y,z,t)$ and $(x',y',z',t')$ respectively. Select $\gamma$ such that $\theta(y+it+\gamma,i) = 0$. Now select $\delta$ such that all other factors in definition of function $g$ don't vanish at the point $v$:
$$
\theta(y'+it'+\delta,i)\theta(y'+it'-\gamma-\delta,i)\neq 0.
$$
It is possible because zeros of theta function are isolated. Also by selecting $\alpha,\beta$ we can assure that function $f$ doesn't vanish at point $v$.

Unless $\theta(y'+it'+\gamma, i)= 0$ the  constructed section solves the  problem. Assume contrary. Since classical theta function has single (up to multiplicity) zero in the fundamental region of lattice, it follows that $y = y'$, $ t = t'$ modulo $\Gamma$.

Select $\alpha$ such that $\theta(z +ix + \alpha,y+i)
= 0$. Note that $\theta(z' +ix' + \alpha,y'+i)
\neq 0$, because otherwise $u=v$. Select $\beta$ such that $f(v)\neq 0$ and $\gamma,\delta$ such that $g(v)\neq
0$.

Thus we constructed necessary section and proved the injectivity of $\varphi_k$.

Let's prove that rank of $\varphi_k$ is maximal. Here we will follow the proof of Lefshetz theorem stated in \cite{Taymanov}.
Firstly we show that the rank of $\varphi_k$ is maximal if and only if the rank of matrix $J$ is maximal
\begin{equation*}
J =
\begin{pmatrix}
  s_1 & & \ldots & & s_{k^2} \\
  \partial_xs_1 & & \ldots & & \partial_xs_{k^2} \\
  \partial_ys_1 & & \ldots & & \partial_ys_{k^2} \\
  \partial_zs_1 & & \ldots & & \partial_zs_{k^2} \\
  \partial_ts_1 & & \ldots & & \partial_ts_{k^2} \\
\end{pmatrix}.
\end{equation*}
The map $\varphi_k$ written in homogenous coordinates is a
composition of the map $\widetilde{\varphi}_k$ to $\mathbb{C}^{k^2}$ and
subsequent projection $\pi: \mathbb{C}^{k^2}\backslash
\{0\}\rightarrow\mathbb{CP}^{k^2-1}$. Obviously if we cross out
first row of $J$ we get the differential of $\widetilde{\varphi}_k$.

Now assume that at the point $u\in {M_{KT}}$ first row of $J$ is a linear combination of other rows.
It means that radius-vector of $\widetilde{\varphi}_k(u)$ is collinear to the image of a certain tangent
vector (to $\M$) at the point $u$. Since $\pi$ projects along complex lines passing through the origin, the kernel of
differential of $\pi$ consists exactly of such vectors.  Therefore rank of $J$ is maximal if and only if rank of $\varphi_k$ is maximal.

Let's transform matrix $J$ to more suitable for us form. The rank of the following matrix coincides with the rank of $J$
\begin{equation*}
\widetilde{J} =
\begin{pmatrix}
  s_1 & & \ldots & & s_{k^2} \\
  (\partial_y-i\partial_t)s_1 & & \ldots & & (\partial_y-i\partial_t)s_{k^2} \\
  (\partial_z-i\partial_x)s_1 & & \ldots & & (\partial_z-i\partial_x)s_{k^2} \\
  (\partial_y+i\partial_t)s_1 & & \ldots & & (\partial_y+i\partial_t)s_{k^2} \\
  (\partial_z+i\partial_x)s_1 & & \ldots & & (\partial_z+i\partial_x)s_{k^2} \\
\end{pmatrix}.
\end{equation*}
Last two rows of $\widetilde{J}$ are the Cauchy-Riemann conditions.
Since sections $s_j$ are holomorphic functions of $z+ix$, the last
row of $\widetilde{J}$ vanishes.

Assume that rank of $\widetilde{J}$ (over $\C$) is less than 4 at the certain fixed
point $u^*=(x^*,y^*,z^*,t^*)\in M_{KT}$. It means that there exist non-trivial constants $a,b,c,d$ such that
\begin{multline*}
as_j(u^*) + \frac{b}{2}(\partial_y-i\partial_t)s_j(u^*) + \frac{c}{2}(\partial_z-i\partial_x)s_j(u^*) +
\frac{d}{2}(\partial_y+i\partial_t)s_j(u^*) = 0,\\\quad j = 1,\ldots ,k^2.
\end{multline*}
By \eqref{product}, the function
$$
s(u,\mu,\nu) = (\mu\cdot\theta_{KT})(u)(\nu\cdot\theta_{KT})(u)((-\mu-\nu)\cdot\theta_{KT})(u).
$$
lies in $\mathcal{L}_k(k=3)$ for any $\mu$,$\nu$. Hence this function is a linear combination of $s_j$ and the identity \begin{multline}
\label{s}
as(u^*,\mu,\nu) + \frac{b}{2}(\partial_y-i\partial_t)s(u^*,\mu,\nu) + \frac{c}{2}(\partial_z-i\partial_x)s(u^*,\mu,\nu) + \\
+ \frac{d}{2}(\partial_y+i\partial_t)s(u^*,\mu,\nu) = 0.
\end{multline}
holds. We define the linear differential operator $L =
\frac{b}{2}(\partial_y-i\partial_t)+\frac{c}{2}(\partial_z-i\partial_x)
+ \frac{d}{2}(\partial_y+i\partial_t)$ and rewrite the last identity
\begin{equation}
\label{a1}
L\log(\mu\cdot\theta_{{KT}})(u^*) = -a - L\log (\nu\cdot
\theta_{{KT}})(u^*)-L\log((-\mu-\nu)\cdot\theta_{{KT}})(u^*).
\end{equation}
For any $u$,$\mu$ there exists $\nu$ such that
\begin{equation}
\label{a2}
(\nu\cdot\theta_{{KT}})(u)((-\mu-\nu)\cdot\theta_{{KT}})(u)\neq 0.
\end{equation}
By \eqref{a1}-\eqref{a2}, the function
\begin{equation}
\label{xi}
\xi(\mu) = L\log(\mu\cdot\theta_{{KT}})(u^*)
\end{equation}
is an entire function of $\mu = (\mu^1,\mu^2)$.
It follows from \eqref{mults} that function $\xi(\mu)$ satisfies the following periodicity conditions
\begin{eqnarray}
\label{localperiod}
\xi(\mu^1 + 1,\mu^2) & = & \xi(\mu^1,\mu^2),\\
\label{localperiod2}
\xi(\mu^1 + y^* + i,\mu^2) & = & \xi(\mu^1,\mu^2) - 2\pi i c,\\
\label{localperiod3}
\xi(\mu^1,\mu^2 + 1) & = & \xi(\mu^1,\mu^2),\\
\label{localperiod4}
\xi(\mu^1 ,\mu^2 + i) & = & \xi(\mu^1,\mu^2) - 2\pi i b.
\end{eqnarray}
Therefore derivatives $\partial_{\mu^j}\xi$ are the periodic and entire functions.
This means that they are constants and $\xi = \alpha\mu^1 + \beta\mu^2 + \gamma$. By \eqref{localperiod},\eqref{localperiod3},
$\alpha=\beta=0$ and function $\xi$ is constant. By \eqref{localperiod2},\eqref{localperiod4}, $b=c=0$. Then the following identity holds
\begin{equation}
\label{d}
\xi(\mu) = \frac{d}{2}\left[\frac{\partial_y\theta(z+ix+\mu^1,y+i)}
{\theta(z+ix+\mu^1, y+i)}\right]_{u=u^*} = \gamma.
\end{equation}
Here we used the Cauchy-Riemann equation implicitly
$$
(\partial_y+i\partial_t)\theta(y+it,i) = 0.
$$
Let's denote by $D$  the  $z+ix$-differentiation
$$
D = \frac{1}{2}(\partial_z-i\partial_x) .
$$
It follows from \eqref{classheat} that
\begin{equation}
\label{dd}
\partial_y\theta(z+ix,y+i) \equiv \frac{1}{4\pi i}(D^2\theta)(z+ix,y+i) - \frac{1}{2}(D\theta)(z+ix,y+i).
\end{equation}
By substituting \eqref{dd} in \eqref{d} and considering that
$$
(D\theta)(z+ix + \mu^1,y+i)= \partial_{\mu^1}\theta(z+ix+\mu^1,y+i),
$$
we have that function $\theta(z^*+ix^*+\mu^1, y^*+i)$ as a function
of $\mu^1$ satisfies the linear ODE with constant coefficients
$$
\frac{d}{4\pi i}\theta'' - \frac{d}{2}\theta' - 2\gamma
\theta = 0.
$$
It is easy to write down general solution of this equation and check that it contradicts  to the periodicity conditions of theta function \eqref{mults}. Thus $d = \gamma = 0$. By \eqref{s}, $a=0$.

As a result all constants $a,b,c,d$ vanish and matrix $\widetilde{J}$ has the maximal rank. Since point $u^*$ is arbitrary, the rank is maximal everywhere.
Theorem proved.

{\sc Remark.} It would be interesting to investigate the connection of this theta function to certain nonlinear equations by the way of  obtaining soliton equations from secant identities for Jacobians \cite{Taymanov}.

\section{Embedding is symplectic}
Manifold $\M$ is a symplectic manifold, where symplectic form  would be for instance the following left-invariant form $\omega_{KT} = (dz-xdy)\wedge dx +
dy\wedge dt$. In this section we will prove the following
\begin{theorem}

\begin{enumerate}
\item If $k\geq 3$ then map $\varphi_k$ induces a symplectic form on $\M$.

\item Induced symplectic form is cohomologous to $k\cdot \omega_{KT}$.
\end{enumerate}
\end{theorem}
{\sc Proof.}\quad Let's choose the following basis of $\mathcal{L}_k$ in the definition of $\varphi_k$
$$
\theta^p_k(z+i x,y+i)\cdot \theta^{q}_k(y+it, i);\quad p,q = 1,\ldots ,k.
$$

Notice that $\varphi_k$ is a composition of two maps.
The first one is $\psi_k:{M_{KT}}\rightarrow\mathbb{CP}^{k-1}\times\mathbb{CP}^{k-1}$, $\psi_k = (\psi_k',\psi_k'')$. Here
$$
\psi_k'(x,y,z) = [\theta^1_k(z+ix,y+i):\ldots:\theta^k_k(z+ix,y+i)],
$$
$$
\psi_k''(y,t) = [\theta^1_k(y+it,i):\ldots:\theta^k_k(y+it,i)].
$$
The second is the Segre map $\sigma_k:\mathbb{CP}^{k-1} \times
\mathbb{CP}^{k-1} \rightarrow \mathbb{CP}^{k^2-1}$, which is defined
in homogenous coordinates by the formula
$$
\sigma_k\left([z^1:\ldots :z^k],[w^1:\ldots :w^k]\right) = [z^1w^1: z^1w^2: \ldots :z^kw^{k-1}: z^kw^k].
$$
Thus $\varphi_k = \sigma_k \circ \psi_k$.
Let's denote by $\Omega'$ symplectic form (associated with the
Fubini--Study metric) on the the first factor of
$\mathbb{CP}^k~\times~\mathbb{CP}^k$, by $\Omega''$ on the second.
Then $\Omega'+\Omega''$ is a symplectic form on the product. Since
the Segre map is a holomorphic embedding, it is sufficient  to show
that induced form $\psi^*_k(\Omega'+ \Omega'')$ is symplectic.

Note that algebra of left-invariant forms on $\M$ is generated by $ dx, dy, dz-xdy, dt $.

Map $\psi_k''$ is holomorphic embedding of complex torus into
$\mathbb{CP}^k$ described by the classical Lefshetz theorem.
Therefore the Fubini--Study form induces symplectic form on torus
$$
{(\psi_k'')^*(y,t)} (\Omega'') = \alpha\cdot dy\wedge dt,
$$
where $\alpha$ doesn't vanish anywhere on $\M$.

Let
$$
{(\psi_k')^*(x,y,z)} (\Omega') = f\cdot (dz-xdy)\wedge dx + g\cdot (dz-xdy)\wedge dy + h\cdot dx\wedge dy.
$$
Here $f,g,h$ are certain functions on $\M$. This is the general view of 2-form on $\M$ generated by map depending on $x,y,z$.

Note that for any fixed $y$ map $\psi_k'$ also becomes holomorphic embedding described by the Lefshetz theorem and therefore
$$
{(\psi_k')}^* (\Omega') =  \beta \cdot dz\wedge dx,
$$
where $\beta$ doesn't vanish anywhere on $\M$.
It follows that $f\equiv \beta$. Gathering altogether
\begin{multline*}
\left(\psi^*_k(\Omega'+\Omega'')\right)^2 = \left({(\psi_k')}^*
(\Omega') + {(\psi_k'')}^*
(\Omega'')\right)^2 =\\
 = \left(\beta\cdot (dz-xdy)\wedge dx + g\cdot (dz-xdy)\wedge dy + h\cdot dx\wedge dy + \alpha\cdot dy\wedge dt \right)^2.
\end{multline*}
By opening the brackets we get that
$$
\left(\psi^*_k(\Omega'+\Omega'')\right)^2 = 2\alpha\beta \cdot
dx\wedge dy\wedge dz \wedge dt.
$$
Last identity means that induced form is non-degenerated. The
closedness is implied by the commutation of differential and
$\psi_k^*$. Thus $\psi^*_k(\Omega'+ \Omega'')$ is closed and
non-degenerate i.e. symplectic form. We proved the first item of the
theorem.

Let's prove the second one. Denote by $L$ the bundle defined by  multiplicators \eqref{mults}.
Earlier we noted that  theta functions of degree $k$ are sections of $L^{\otimes k}$.

Recall that any complex line bundle over a manifold $M$ is induced by
the universal bundle over the complex projective space when $M$ is
mapped to $\mathbb{CP}^n$. Therefore $L^{\otimes k}$ is a pullback
of universal bundle and curvature form of $L^{\otimes k}$ is a
pullback of the Fubini--Study form, which is a curvature form of
universal bundle.  Recall also that the first Chern class of line
bundle is realized by curvature form. Thus the cohomological class
of induced form coincides with $c_1(L^{\otimes k}) = k\cdot c_1(L)$
and we must prove that
$$
c_1(L) = [(dz-xdy)\wedge dx + dy\wedge dt].
$$
Consider the cover of $\M$ by sets
$$
U_{\lambda} = \lambda\cdot U_0,\quad \lambda\in\Gamma
$$
where $U_0 = \{|u^k|<3/4\}$.
The nerve $N(\mathcal{U})$ of the minimal subcover $\mathcal{U}$ of the above cover $U_\lambda$ is homeomorphic to $\M$ and
its cohomologies with coefficients in $\Z$ coincide with $H^*(\M;\Z)$.

The coordinate transformations $g_{\lambda\mu}:U_{\lambda}\cap
U_{\mu}\rightarrow \mathbb{C}^*$ are expressed  in terms of multiplicators
\begin{equation}
\label{g}
g_{\lambda\mu}(u) = e_{\lambda}(u)\cdot e_{\mu^{-1}}(\mu\cdot u);\quad \lambda,\mu\in\Gamma.
\end{equation}
By definition the cocycle
$z_{\lambda\mu\nu} \in C^2(\mathcal{U};\mathbb{Z})$
\begin{equation}
\label{cocycle}
z_{\lambda\mu\nu} = \frac{1}{2\pi i}(\log(g_{\lambda\mu}) + \log(g_{\mu\nu}) - \log(g_{\nu\lambda}))
\end{equation}
realizes the first Chern class of bundle $L$. Given formula is the value of $z$ at
a two-dimensional simplex $(\lambda,\mu,\nu)\in N(\mathcal{U})$.

The group $H_2(\M;\Z)$ is generated by homological classes of tori
$T_{ac}$, $T_{bc}$, $T_{da}$, $T_{db}$, spanned by commuting
translations \eqref{period}.

Define the functions $f_\lambda(u)$ as follows
\begin{equation}
\label{f}
e_\lambda(u) = e^{2\pi i f_\lambda(u)}.
\end{equation}
By \eqref{g}--\eqref{f},
\begin{equation*}
c_1([T_{\lambda\mu}]) = f_\mu(u) + f_\lambda(\mu \cdot u) - f_\lambda(u) - f_\mu(\lambda\cdot u).
\end{equation*}
Calculate the first Chern class at the basis tori
\begin{equation}
\label{c1calc}
c_1([T_{ca}]) = c_1([T_{bd}]) = 1;\  c_1([T_{cb}]) =  c_1([T_{ad}]) = 0.
\end{equation}
Since manifold $\M$ is a homogenous space of nilpotent Lie group, any element of $H^2({M_{KT}};\mathbb{R})$ is realized by the left-invariant form dual to basis cocycle. Group $H^2({M_{KT}};\mathbb{R})$ is generated by cohomological classes of the forms $(dz-xdy)\wedge dx$, $dy\wedge dt$, $(dz-xdy)\wedge dy$ and
$dx\wedge dt$.

By \eqref{c1calc}, $c_1(L) = [(dz-xdy)\wedge dx +
dy\wedge dt]$. Theorem is proved.

\noindent
Institute of Mathematics and Informatics, Yakut State University, Yakutsk, Russia.\\
E-mail:\ {\tt egorov.dima@gmail.com}

\end{document}